\documentclass[a4paper,11pt]{article}
\usepackage{latexsym}
\usepackage{amssymb}
\usepackage{amsfonts}
\usepackage{amsmath}
\usepackage{indentfirst}
\usepackage{xcolor}
\usepackage{xypic}
\usepackage{graphicx}
\usepackage{epsfig}
\newtheorem{prop}{Proposition}[section]
\newtheorem{teor}{Theorem}[section]
\newtheorem{lemma}{Lemma}[section]
\newtheorem{cor}{Corollary}[section]
\newcommand{\cvd}{\hfill $\blacksquare$\bigskip}

\newcommand{\ninN}{n\in \mathbf{N}}

\definecolor{myblue}{rgb}{0,0.1,0.9}

\definecolor{myred}{rgb}{0.9,0.1,0}

\date{}
\author{Filippo Disanto\thanks{Institut
f\"ur Genetik, Universit\"at zu K\"oln, Z\"ulpicher Str. 47a,
50674 K\"oln, Germany \tt{fdisanto@uni-koeln.de}} \and Luca Ferrari\thanks{Universit\'a di Firenze, Dipartimento di
Sistemi e Informatica, viale Morgagni 65, 50134 Firenze, Italy
\tt{ferrari@dsi.unifi.it}}\and Simone Rinaldi\thanks{Universit\'a di
Siena, Dipartimento di Scienze Matematiche ed Informatiche, Pian
dei Mantellini, 44, 53100, Siena, Italy \tt{rinaldi@unisi.it}}}
\title{A partial order structure on interval orders}

\begin{document}

\maketitle

\begin{abstract}
We introduce a partial order structure on the set of interval
orders of a given size, and prove that such a structure is in fact
a lattice. We also provide a way to compute meet and join inside
this lattice. Finally, we show that, if we restrict to series
parallel interval order, what we obtain is the classical Tamari
poset.
\end{abstract}

\bigskip

%

\section{Introduction}

\emph{Interval orders} are an interesting class of partial orders,
introduced by Fishburn in \cite{F1}, which are especially
important even in non strictly mathematical contexts, such as
experimental psychology, economic theory, philosophical ontology
and computer science \cite{F2}. From a purely combinatorial point
of view, some remarkable features of interval orders have been
recently exploited in \cite{BMCDK}, where their connection with
some interesting combinatorial structures, such as pattern
avoiding permutations and chord diagrams, have been shown.
Starting from that paper, a number of articles has been written,
trying to go deeper in the combinatorial knowledge of interval
orders.

In our work we will explore the possibility of introducing a
suitable partial order structure on the set of interval orders
(having ground set of fixed size). Our goal is twofold: the
resulting poset should be as ``natural" as possible, and it should
be compatible with possible (already existing) partial orders on
subsets of its ground set. We have been able to fulfill this goal,
by defining a presumably new partial order structure which is
easily defined in terms of a very natural labelling of the
elements of the ground set, which has the additional features of
being a lattice and of coinciding with the Tamari order when
restricted to series parallel interval orders.

The article is organized as follows. In section \ref{notations} we
recall those definitions and facts concerning interval orders and
poset theory in general that will be useful throughout the paper.
In section \ref{label} we introduce a particular labelling of an
interval order that will be crucial for the definition of our
partial order on interval orders of the same size. Section
\ref{sectionalgebr} is the heart of the paper, and contains the
proof that our poset is in fact a lattice. Finally, section
\ref{tamari} provides the argument to show that our partial order,
restricted to series parallel interval orders, is isomorphic to
the well-known Tamari order.

\section{Notations and preliminaries}\label{notations}

Let $P=(X,\leq )$ be a finite poset. A \emph{linear extension of}
$P$ is a bijection $\lambda :X \rightarrow \{ 1,2,\dots ,|X|\}$
such that $x<y$ in $P$ implies $\lambda (x)<\lambda (y)$.

Given $Y\subseteq X$, the \emph{up-set of} $P$ \emph{generated by}
$Y$ is the set $F_P (Y)=\{ x\in X\; |\; \forall y\in Y,x>y\}$.
Analogously, the \emph{down-set of} $P$ \emph{generated by} $Y$ is
the set $I_P (Y)=\{ x\in X\; |\; \forall y\in Y,x<y\}$. In
particular, we will denote with $\mathcal{D}_P =\{ I_P (\{ x\} )\;
|\;x\in X\}$ and $\mathcal{U}_P =\{ F_P (\{ x\} )\; |\; x\in X\}$
the sets of \emph{principal down-sets} and \emph{principal
up-sets} of $P$, respectively. To simplify notations, we will
often write $I(x)$ in place of $I_P (\{ x\} )$ and, analogously,
$F(x)$ in place of $F_P (\{ x\} )$.

Observe that the above definitions of an up-set and of a down-set
slightly differ from the usual ones which can be found in the
literature. Indeed, in this work an up-set generated by $Y$ does
not contains the elements of $Y$ (and the same convention holds
for down-sets). We have preferred to give definitions in this way
since this will help us in stating (and then proving) our main
results.

Given $x,y\in X$, we say that $x$ and $y$ are \emph{order
equivalent} whenever $I(x)=I(y)$ and $F(x)=F(y)$. In this case, we
will use the notation $x\sim y$.

We say that a poset $P$ \emph{avoids} a poset $S$ when $P$ has no
subposet isomorphic to $S$. Borrowing notations from the theory of
pattern avoiding permutations, we will refer to the class of
posets avoiding the poset $S$ using the symbol $AV(S)$; in
particular, when we restrict ourselves to posets of cardinality
$n$, we will write $AV_n(S)$.

%
%

An important class of posets is that of interval orders
\cite{BMCDK,EZ,F1,Kh}. A poset $P=(X,\leq )$ is called an
\emph{interval order} when there exists a function $J$ mapping
each element $x\in X$ into a closed interval $J(x)=[a_x ,b_x
]\subseteq \mathbf{R}$ in such a way that, for all $x,y\in X$,
$x<y$ in $P$ if and only if $b_x <a_y$ in $\mathbf{R}$. We call
$J$ an \emph{(interval) representation} of $P$. If the interval
order $P$ is finite, then we can obviously find a representation
of $P$ such that, for every element $x$, the values $a_x$ and
$b_x$ are integers.

\begin{figure}[ht]
\begin{center}
\centerline{\hbox{\psfig{figure=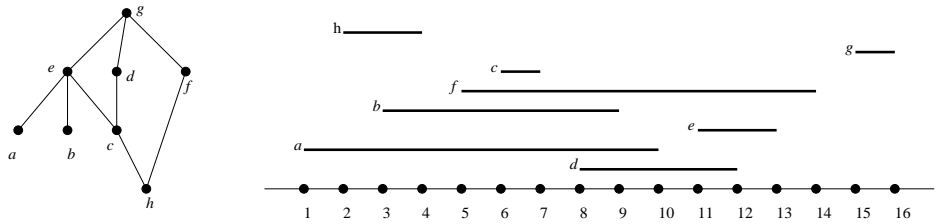,width=4.8in,clip=}}}
\caption{An interval order and one of its representations.}
\label{representation}
\end{center}
\end{figure}

\begin{figure}[htb]
\begin{center}
\centerline{\hbox{\psfig{figure=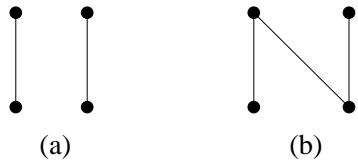,width=1.8in,clip=}}}
\caption{(a) The poset $2+2$; (b) The fence of order four.}
\label{poset}
\end{center}
\end{figure}

In \cite{F1} Fishburn gives the following characterization for the
class of interval orders in terms of avoiding subposets. Recall
that the poset $2+2$ is the disjoint union on two chains each
having two elements (see Figure \ref{poset}).

\begin{teor}\label{fish}
A poset $P = (X,\leq)$ is an interval order if and only if $P\in
AV(2+2)$.
\end{teor}

The following proposition, stated in \cite{Kh}, gives a
characterization for the class of interval orders in terms of
principal down-sets and principal up-sets.

\begin{prop}\label{eq'}
The following statements are equivalent:
\begin{enumerate}
\item[i)] $P$ is an interval order;

\item[ii)] any two distinct sets in $\mathcal{D}_P$ are ordered by
inclusion;

\item[iii)] any two distinct sets in $\mathcal{U}_P$ are ordered
by inclusion.
\end{enumerate}
\end{prop}

\section{The admissible labelling of an interval order}\label{label}

Let $P(X,\leq )$ be a poset. The following proposition gives an
immediate characterization of order equivalent elements, whose
easy proof is left to the reader.

\begin{prop}\label{auto} Two elements $x$ and $y$ of a poset $P$ are order
equivalent if and only if the map from $P$ to itself which
exchanges $x$ and $y$ is an automorphism of $P$.
\end{prop}

From now on in this section, the poset $P$ will denote an interval
order.

\bigskip

A linear extension $\lambda$ of $P$ is called an \emph{admissible
labelling} of $P$ whenever, for all $x,y\in X$, if $\lambda
(x)<\lambda (y)$ then either $I(x)\subset I(y)$, or $I(x)=I(y)$
and $F(x)\subset F(y)$, or $x\sim y$.
%
%

Such a labelling has been defined in \cite{DR}, 
in the context of a recursive construction of interval orders,
where it is also shown that each interval order admits at least
one admissible labelling.

A trivial property of an admissible labelling of an interval order
(which will be useful in the next section) is the following.

\begin{prop}\label{trivial} Let $\lambda$ be an admissible labelling of $P$.
Given $x,y,z\in X$ such that $x<y$ and $\lambda (y)\leq \lambda
(z)$, then $x<z$.
\end{prop}

\emph{Proof.}\quad From the definition of admissible labelling,
$\lambda (y)\leq \lambda (z)$ implies that $I(y)\subseteq I(z)$.
Since $x\in I(y)$, we have that $x\in I(z)$, that is $x<z$.\cvd

\bigskip

The rest of this section is devoted to show that an interval order
admits a unique admissible labelling (up to automorphisms).

\begin{lemma}\label{<>}
Suppose that $\lambda_1$ and $\lambda_2$ are two admissible
labellings of $P$. If $\lambda_1(x)>\lambda_1(y)$ and
$\lambda_2(x)<\lambda_2(y)$, then $x \sim y$.
\end{lemma}
\emph{Proof.}\quad This follows immediately from the definition of
an admissible labelling. \cvd

\begin{prop} \label{lambde}
Suppose that $\lambda_1$ and $\lambda_2$ are two admissible
labellings of $P$. If $\lambda_1(x)=\lambda_2(y)$, then $x \sim
y$.
\end{prop}

\emph{Proof.}\quad Let $\ell=\lambda_1 (x)=\lambda_2 (y)$.
Concerning the values of the two labels $\lambda_2 (x)$ and
$\lambda_1 (y)$ we have essentially two different cases.

\begin{enumerate}

\item If $\lambda_2 (x), \lambda_1 (y)<\ell$, then we can simply
apply the above lemma. The same argument can be used in the case
$\lambda_2 (x), \lambda_1 (y)>\ell$.

\item Suppose, without loss of generality, that $\lambda_1
(y)<\ell <\lambda_2 (x)$. We then claim that there exists $z\in X$
such that $z\sim x$ and $z\sim y$ (whence the thesis will easily
follow by transitivity). Indeed, we observe that
\begin{displaymath}
|\{ z\in X\setminus \{ x,y\} \; |\; \lambda_2 (z)>\ell \} |<|\{
z\in X\setminus \{ x,y\} \; |\; \lambda_1 (z)>\ell \} |
\end{displaymath}
(since, in the labelling $\lambda_2$, the label of $x$ is greater
than $\ell$). Thus there exists an element $z\in X\setminus \{
x,y\}$ such that $\lambda_1 (z)>\ell$ and $\lambda_2 (z)<\ell$.
Since we are supposing that $\ell <\lambda_2 (x)$, we then have
that $\lambda_1 (z)>\lambda_1 (x)$ and $\lambda_2 (z)<\lambda_2
(x)$. Therefore we can apply once again lemma \ref{<>} to obtain
that $z\sim x$. An analogous argument shows that $z\sim y$.\cvd

\end{enumerate}


\begin{cor} \label{unica}
Let $\lambda_1 ,\lambda_2$ be two admissible labellings of $P$.
Then there exists an automorphism $f$ of $P$ such that, for all
$x\in P$, $\lambda_1 (x)=\lambda_2 (f(x))$.
\end{cor}

\emph{Proof.}\quad Given $x\in X$, let $f(x)$ be the (unique)
element $y\in X$ such that $\lambda_1 (x)=\lambda_2 (y)$. Thanks
to proposition \ref{lambde}, we have $x\sim y$, and so the map $f$
is an automorphism of $P$ (since it is the composition of
automorphisms, by proposition \ref{auto}).\cvd

%
%
%
%
%
The last corollary states that there exists a unique admissible
labelling of a given interval order \emph{up to order
automorphism}. This uniqueness result will be frequently used in
the rest of the paper.

\section{The poset $AV_n(2+2)$}\label{sectionalgebr}

In the present section, which is the heart of our work, we endow
each set $AV_n(2+2)$ with a partial order structure. We then prove
that the resulting poset is in fact a lattice, which provides a
generalization of the \emph{Tamari lattice}. This partial order on
$AV_n(2+2)$ is believed to be new.

\bigskip

In the sequel we will consider \emph{interval orders endowed with
their admissible labelling}, and we will identify an element $x$
with its label
$\lambda (x)$ (this can be done thanks to Corollary \ref{unica}). Moreover, 
for an interval order $P=(X,\leq )$, we will write $x\leq y$ to
indicate that the element $x$ is less than or equal to $y$ with
respect to the partial order of $P$, whereas we will write $x\prec
y$ to mean that the label of $x$ is less than the label of $y$.
See Figure \ref{equivalence1} for an example.

\bigskip


Given two interval orders $P_1$ and $P_2$ on the same ground set
$X=[n]=\{ 1,2,\dots ,n\}$, we declare $P_1 \leq_T P_2$ whenever
$P_1 \supseteq P_2$, i.e. the (partial order) relation $P_2$ is a
subset of the (partial order) relation $P_1$.


The following proposition characterizes the order relation
$\leq_T$ in terms of both the principal up-sets and the principal
down-sets of the elements of $AV_n(2+2)$. The proof is an easy
consequence of the notations and results previously recalled, so
it is left to the reader.

\begin{prop}\label{rd1}
Let $P_1 ,P_2$ be two interval orders on $X=\{1,2,\dots,n \}$. The
following are equivalent:
\begin{enumerate}
\item[i)] for each $x \in X$, $F_{P_1}(x)\supseteq F_{P_2}(x)$;
\item[ii)] for each $x \in X$, $I_{P_1}(x)\supseteq I_{P_2}(x)$;
\item[iii)]$P_1 \leq_T P_2$
\end{enumerate}
\end{prop}

Let $\Gamma =\{ (i,j)\in [n]\times [n]\; |\; i\leq j\}$. Then the
set of interval orders on $[n]$ is clearly what is usually called
a \emph{family of subsets of} $\Gamma$. We recall here a classical
definition which can be found, for instance, in \cite{DP}. A
family of subsets $\mathcal{L}$ of a set $\Gamma$ is called a
\emph{closure system} on $\Gamma$ when it is closed under
arbitrary intersections and it contains $\Gamma$. Analogously,
when $\mathcal{L}$ is closed under arbitrary unions and it
contains the empty set, it will be called a \emph{dual closure
system} on $\Gamma$.

The following result (recorded in \cite{DP} as well) gives an
important property of closure systems.

\begin{teor} Any closure system
$\mathcal{L}$ 
is a complete lattice, in which
\begin{eqnarray}\label{topped}
\bigwedge_{i\in I}A_i &=&\bigcap_{i\in I}A_i ,\nonumber \\
\bigvee_{i\in I}A_i &=&\bigcap\{ B\in \mathcal{L}\; |\; B\supseteq
\bigcup_{i\in I}A_i\},
\end{eqnarray}
for all $\{ A_i \}_{i\in I}\subseteq \mathcal{L}$.

Analogously, any dual closure system $\mathcal{L}$ 
is a complete lattice, in which
\begin{eqnarray}\label{dualtopped}
\bigwedge_{i\in I}A_i &=&\bigcup_{i\in I}A_i ,\nonumber \\
\bigvee_{i\in I}A_i &=&\bigcup\{ B\in \mathcal{L}\; |\; B\subseteq
\bigcap_{i\in I}A_i\},
\end{eqnarray}
for all $\{ A_i \}_{i\in I}\subseteq \mathcal{L}$.
\end{teor}

In view of the above theorem, the following result is trivial, so
it is stated without proof.

\begin{lemma}\label{prel} Let $\mathcal{L}$ be a family of subsets of a set
$\Gamma$. Suppose that there exists $A\in \mathcal{L}$ such that
$A\subseteq B$ for all $B\in \mathcal{L}$ and $\mathcal{L}$ is
closed under arbitrary nonempty unions. Then $\mathcal{L}$ is a
complete lattice, in which the meet and join operations are
computed as in (\ref{dualtopped}).
\end{lemma}

The above facts allow us to formulate our main result concerning
the order structure of $AV_n(2+2)$.

\begin{teor} $(AV_n(2+2),\leq_T)$ is a (complete) lattice, in which the meet and join operations are expressed as follows:
\begin{eqnarray*}
P_1 \wedge P_2 &=&P_1 \cup P_2 ,\\
P_1 \vee P_2 &=&\bigcup \{ P\in AV_n (2+2)\; |\; P \subseteq P_1 \cap P_2 \}.
\end{eqnarray*}
\end{teor}

\emph{Proof.}\quad We start by observing that $D=\{ (i,i)\; |\;
i\in [n]\}$ is an interval order on $[n]$ (since it is the
discrete poset on $[n]$) and that any interval order on $[n]$
clearly contains $D$. Thus $AV_n(2+2)$ is a family of subsets of
$\Gamma$ having $D$ as a minimum. Therefore, since $AV_n(2+2)$ is
finite, in view of Lemma \ref{prel} it will be enough to prove
that, for any $P_1 ,P_2 \in AV_n(2+2)$, $P=P_1 \cup P_2 \in
AV_n(2+2)$. In what follows, we will denote by $\leq_P$ the
partial order relation on $P$, and by $\leq_{P_i }$ the partial
order relation on each $P_i$, for $i=1,2$.

The first thing to prove is that $P$ is a poset. In fact, $P$ is
trivially reflexive (since it contains $D$). Moreover, suppose
that $x\leq_P y$ and $y\leq_P x$. If the two relations hold in the
same poset $P_i$ (that is, if $x\leq_{P_i }y$ and $y\leq_{P_i }x$
for $i=1$ or $i=2$), then trivially $x=y$. Otherwise, suppose
(without loss of generality) that $x\leq_{P_1 }y$ and $y\leq_{P_2
}x$. Since the admissible labelling is a linear extension of its
interval order, then necessarily $x\preceq y$ and $y\preceq x$,
whence immediately $x=y$, and so $\leq_P$ is antisymmetric.
Finally, suppose that $x\leq_P y$ and $y\leq_P z$. Also in this
case, the only nontrivial case arises when (without loss of
generality) $x<_{P_1}y$ and $y<_{P_2}z$. In particular, this
implies that $y\preceq z$. Thus, thanks to proposition
\ref{trivial}, we can conclude that $x<_{P_1}z$, whence $x<_P z$,
that is $\leq_P$ is transitive.

Our next goal is to show that $P$ is an interval order. Thanks to
proposition \ref{eq'}, we will achieve this by showing that, if
$x\preceq y$, then $I_P (x)\subseteq I_P (y)$. Indeed, let $z\in
I_P (x)$, i.e. $z<_P x$. Without loss of generality, this means
that $z<_{P_1 }x$. Together with $x\preceq y$, thanks to
proposition \ref{trivial}, this implies that $z<_{P_1 }y$, and so
$z<_P y$, i.e. $z\in I_P (y)$.

Finally, we observe that the labelling of the elements of $P$
induced by $P_1$ and $P_2$ is an admissible labelling. Indeed, it
is easy to show (and so left to the reader) that such a labelling
is a linear extension of $P$, and that, for each $x$, $I_P
(x)=I_{P_1}(x)\cup I_{P_2}(x)$ and $F_P (x)=F_{P_1}(x)\cup
F_{P_2}(x)$.\cvd

\begin{figure}
\begin{center}
\centerline{\hbox{\psfig{figure=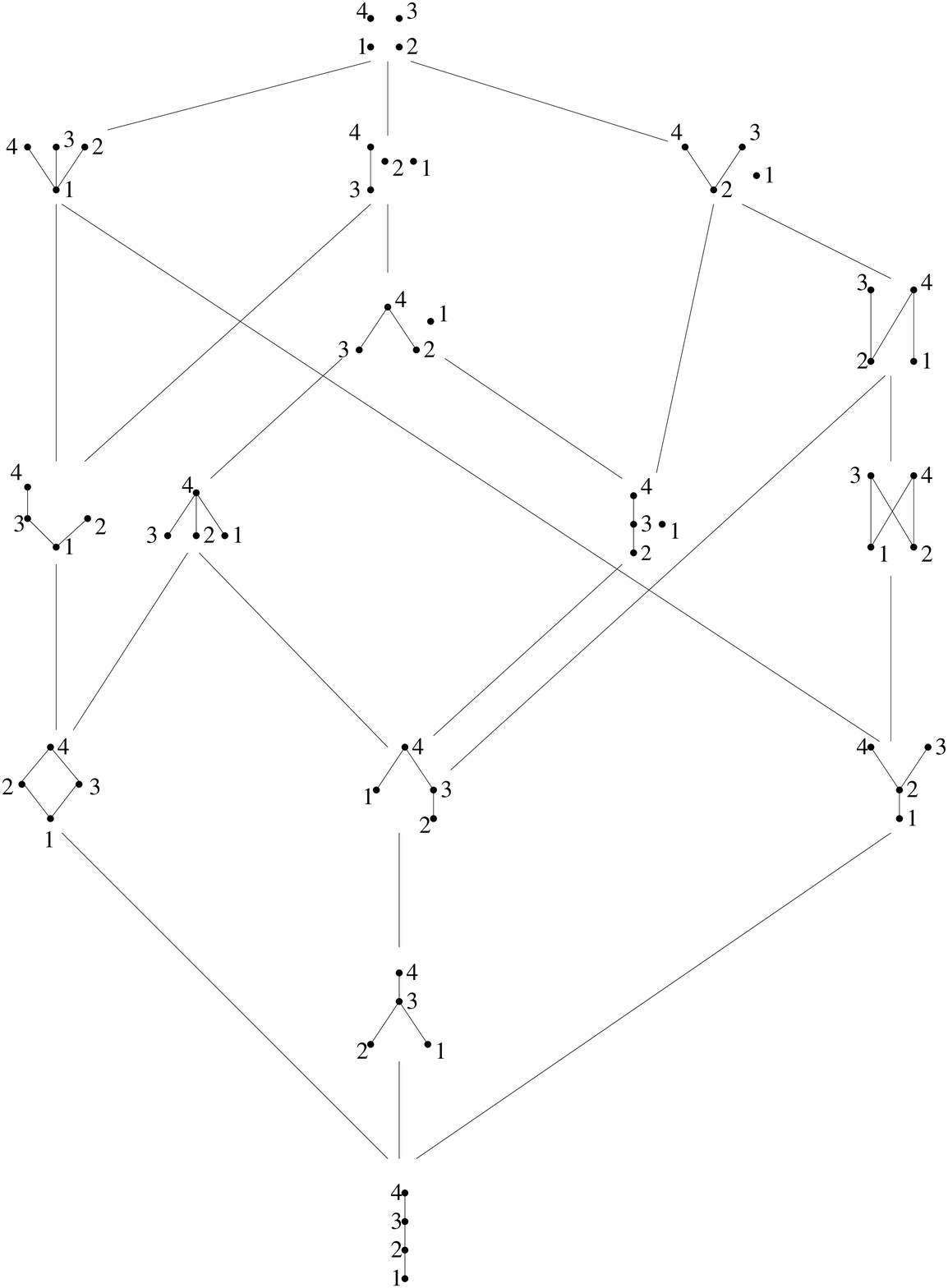,width=3.2in,clip=}}}
\caption{The Hasse diagram of the lattice $(AV_4 (2+2),\leq_T )$}
\label{5}
\end{center}
\end{figure}

\section{The Tamari lattice on series parallel interval
orders}\label{tamari}

%
%
%
%
%

In this final section we will consider the restriction of the
poset $(AV_n (2+2),\leq_T )$ to the set of series parallel
interval orders. This means that we will focus on the poset $(AV_n
(2+2,N),\leq_T )$, where $N$ denotes the \emph{fence having four
elements} (see Figure~\ref{poset}). In particular we will show
that, for any positive integer $n$, $(AV_n (2+2,N),\leq_T )$ is
the Tamari lattice of order $n$.

We point out that, in the literature, there are several extensions
of the Tamari lattice, see for instance \cite{R,S,T}. However, to
the best of our knowledge, the extension we propose in this paper
does not match any of them.

\bigskip

\begin{figure}[ht]
\begin{center}
\centerline{\hbox{\psfig{figure=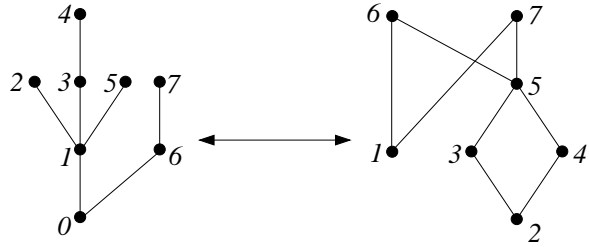,width=3.0in,clip=}}}
\caption{A planar tree with its preorder labelling and the
associated series parallel interval order with its admissible
labelling.} \label{equivalence1}
\end{center}
\end{figure}

From now on, we will consider planar rooted trees whose nodes are
labelled according to the preorder visit (with the root labelled
$0$) and we will systematically identify a node of a tree with its
label (as in the tree represented in Figure \ref{equivalence1}).
Moreover, we will write $x\prec y$ to mean that the label of the
node $x$ is less than the label of the node $y$. Finally, we will
depict trees with their root at the bottom; so, words like
\emph{left} or \emph{right} will refer to this representation (in
particular, the sons of a node will be canonically ordered from
left to right). Given a planar tree $T$ and one of its nodes $k$,
let $u_T (k)$ 
be 
the set of
\emph{descendants} 
of $k$ in the tree $T$. 


Let $T$ be a planar tree and $[n]=\{ 1,2,\ldots ,n\}$ the set of
its nodes different from the root. We define a binary relation $R$
on $[n]$ by setting $xRy$ whenever either $x=y$ or the following
two facts hold: $y\notin u(x)$ and $x\prec y$. This map has been
considered in \cite{DFPR} to define a lattice structure on series
parallel interval order which is isomorphic to the Tamari lattice,
and is illustrated in Figure \ref{equivalence1}. In particular, in
the above cited paper, the authors proved a bunch of results that
are collected in the following proposition.


\begin{prop} \label{ricordi}
\begin{itemize}

\item[$1.$] The structure $P_T =([n],R)$ is a series parallel
interval order, and the labelling of its elements is a linear
extension of $R$.

\item[$2.$] Every series parallel interval order of size $n$ is
isomorphic to $P_T$, for some planar tree $T$.

\item[$3.$] Set $P_1=P_{T_{1}}$ and $P_2=P_{T_{2}}$, if we define
$P_1 \leq_t P_2$ when, $\forall x \in [n]$, $I_{P_1}(x) \supseteq
I_{P_2}(x)$, then $(AV_n(2+2,N),\leq_t)$ is the Tamari lattice of
order $n$.

\end{itemize}
\end{prop}

Our goal is now to show that, if we restrict the order relation
$\leq_T$ defined in the previous section to the set of series
parallel interval orders,
we obtain precisely the Tamari poset.

\begin{prop} The labelling of the poset $P_T$ determined by the preorder
visit on the associated tree $T$ coincides with the admissible
labelling of $R$.
\end{prop}

\emph{Proof.}\quad We start by observing that the labelling of
$P_T$ determined by the preorder visit on $T$ is indeed a linear
extension of $R$ (thanks to Proposition \ref{ricordi}), so the
statement of this proposition makes sense.

Consider $x,y\in [n]$, with $x\preceq y$. We first observe that,
in $P_T$, $I(x)\subseteq I(y)$, that is, for all $z\in [n]$, $zRx$
implies $zRy$ (the proof of this assertion is very easy, so we
leave it to the reader). Now suppose that $x\prec y$ and
$I(x)=I(y)$. This implies that $y\in u(x)$, since otherwise we
would have $xRy$, i.e. $x\in I(y)$, which is not possible (recall
that $x\notin I(x)$). If $xRz$, then necessarily $z\notin u(x)$,
and so a fortiori $z\notin u(y)$; moreover, it is immediate to see
that $y\prec z$. Therefore we have $yRz$, thus proving that
$F(x)\subseteq F(y)$.\cvd

The above proposition shows that the notion of admissible
labelling, when restricted to series parallel interval orders,
coincides with the notion of \emph{preorder linear extension}
introduced in \cite{DFPR}. Thus, using  $3$ of
Proposition~\ref{ricordi}, we can finally state the main theorem
of this section.

\begin{teor}\label{ref3}
The Tamari lattice of order $n$ is the restriction of the lattice
$(AV_n (2+2),\leq_T )$ to the set of series parallel interval
orders $AV_n (2+2,N)$.
\end{teor}

\section{Further work}

The main aim of the present work has been the definition of a
(presumably new) lattice structure on interval orders which,
restricted to series parallel interval orders, turns out to be
isomorphic to the classical Tamari lattice structure. However,
concerning the general (order-theoretic) properties of such a
structure, we have only scratched the surface, and we believe that
it would be very interesting to go deeper into the knowledge of
these lattices. As an example of what could be done, we close our
paper with a structural result which gives some insight on the
relationship between the poset of interval orders and its subposet
of series parallel interval orders.

\begin{prop} For every $\ninN$, $AV_n (2+2,N)$ is a
meet subsemilattice (but not in general a join subsemilattice) of
$AV_n (2+2)$.
\end{prop}

\emph{Proof.}\quad The fact that $AV_n (2+2,N)$ is not in general
a join subsemilattice of $AV_n (2+2)$ can be easily verified, for
instance, by inspecting Figure \ref{5} (and by noticing that the
fence of order four can be obtained as the join of two series
parallel interval orders).

In order to prove that $AV_n (2+2,N)$ is a meet subsemilattice of
$AV_n (2+2)$, we argue by contradiction and suppose that, given
$P=P_1 \wedge P_2$, with $P_1 ,P_2 \in AV_n (2+2,N)$, there exists
a subposet of $P$ isomorhpic to $N$ ($P$ cannot contain any
subposet isomorphic to $2+2$, of course). To fix notations,
suppose that $\{ a,b,c,d\}$ is an occurrence of the poset $N$
inside $P$, with $a\leq_P c,b\leq_P c$ and $b\leq_P d$. It is
clear that there cannot exist $i\in \{1,2\}$ such that the above
listed inequalities hold in $P_i$. Thus we have three essentially
distinct cases.

\begin{itemize}

\item[a)] $a\leq_{P_1}c,b\leq_{P_2}c$ and $b\leq_{P_1}d$. This
case is plainly impossible, otherwise $\{a,b,c,d\}$ would be an
occurrence of $2+2$ inside $P_1$.

\item[b)] $a\leq_{P_1}c,b\leq_{P_1}c$ and $b\leq_{P_2}d$. In this
case, in $P_1$ $d$ is incomparable with any of the remaining three
elements (otherwise there would be an occurrence of $N$ in $P_1$).
We claim that, in $P_1$, $I(d)\subseteq I(b)$. Indeed, if we had
$x\in I(d)$ and $x\notin I(b)$, then there would exist $x$ such
that $x<d$ and $x\nless b$. A simple argument then shows that $\{
b,c,x,d\}$ would constitute an occurrence of either $N$ or $2+2$
in $P_1$ (depending on whether $x<c$ or not), which is not
possible. From $I(d)\subseteq I(b)$ and $F(b) \nsubseteq F(d)$ 
we deduce that $d\prec b$, but this leads to a contradiction,
since we are supposing that $b\leq_{P_2}d$ (which implies that $b\prec d$).

\item[c)] $a\leq_{P_1}c,b\leq_{P_2}c$ and $b\leq_{P_2}d$. In this
case, we can assume that both $b$ and $d$ are incomparable with
any of the three remaining elements in $P_1$ (otherwise one of the
above cases would occur). We claim that, in $P_1$, $I(b)\subseteq
I(a)$. Indeed, if we had $x\in I(b)$ and $x\notin I(a)$, then
there would exist $x$ such that $x<b$ and $x\nless a$. A simple
argument then shows that $\{ a,b,c,x\}$ would constitute an
occurrence of either $N$ or $2+2$ in $P_1$ (depending on whether
$x<c$ or not), which is not possible. From $I(b)\subseteq I(a)$ we
deduce that $b\prec a$. Moreover, in $P_2$, $a$ is easily seen to
be incomparable with the remaining three elements; from this fact,
using an argument similar to the previous ones (and whose details
are then left to the reader), we deduce that $I(a)\subseteq I(b)$,
whence we get $a\prec b$ since also $F(b) \nsubseteq F(a)$, which
contradicts what previously shown.

\end{itemize}

Thus we have shown that, in all cases, $P$ cannot contain any
occurrence of the subposet $N$, which was enough to conclude.\cvd

\end{document}